*Correction*

# Correction to: Equidistribution of Eisenstein Series in the Level Aspect

**Ikuya Kaneko**[1], **Shin-ya Koyama**[2]

[1] Tsukuba Kaisei High School, 3315-10 Kashiwada, Ushiku, Ibaraki 300-1211, Japan.
E-mail: ikuyak@icloud.com
[2] Department of Biomedical Engineering, Toyo University, 2100 Kujirai, Kawagoe, Saitama 350-8585, Japan. E-mail: koyama@tmtv.ne.jp

**Correction to: Commun. Math. Phys. https://doi.org/10.1007/s00220-009-0764-x**

The second author formulated quantum unique ergodicity for Eisenstein series in the prime level aspect in “Equidistribution of Eisenstein series in the level aspect”, Commun. Math. Phys. **289**(3), 1131–1150 (2009). We point out errors and correct the proofs with partially weakened claims.

This article of corrigendum and addendum comprises three sections. Section 1 is devoted to making concise corrections to the claim and the proof of Theorems 1.2 and 1.3, and Sect. 2 serves as a modification of the computation of the contribution from Maaß cusp forms. Other relatively marginal mathematical issues and notational/typographical errors are listed afterwards in Sect. 3.

## 1. Correction to Theorem 1.3

In this section we correct the main theorem (Theorem 1.3). In the original paper [7], Theorem 1.2 was derived as a lucid corollary of Theorem 1.3, while we conclude that the content of Theorem 1.2 ought to hold with a little modification. In our fresh proof, we could no longer rely upon the Luo–Sarnak approach as the level is varying. It behoves us to mention that the quantum unique ergodicity (QUE) conjecture for unitary Eisenstein series in the level aspect actually *cannot* hold for all systems of test functions whose levels are sufficiently large. This is in contradiction to the results shown in [7], where the test function was supposed to be of level $q$.

We define Eisenstein series indexed by cusps $\mathfrak{a}$ singular for a nebentypus $\chi$:

$$E_{\mathfrak{a}}^{(q)}(z, s, \chi) = E_{\mathfrak{a}}(z, s, \chi) := \sum_{\gamma \in \Gamma_{\mathfrak{a}} \backslash \Gamma_0(q)} \overline{\chi}(\gamma) \Im(\sigma_{\mathfrak{a}}^{-1} \gamma z)^s,$$

The original article can be found online at https://doi.org/10.1007/s00220-009-0764-x.

which converges absolutely for $\Re(s) > 1$ and $z \in \mathbb{H}$, where $\Gamma_{\mathfrak{a}} := \{\gamma \in \Gamma_0(q) : \gamma\mathfrak{a} = \mathfrak{a}\}$ stands for the stabiliser of the cusp $\mathfrak{a}$, and the scaling matrix $\sigma_{\mathfrak{a}} \in \mathrm{SL}_2(\mathbb{R})$ is such that $\sigma_{\mathfrak{a}}\infty = \mathfrak{a}$, $\sigma_{\mathfrak{a}}^{-1}\Gamma_{\mathfrak{a}}\sigma_{\mathfrak{a}} = \Gamma_\infty$. In the following lines, we often drop the superscript $(q)$, and also suppress $\chi$ if it is principal. So long as there is no risk of confusion, we simply write the Eisenstein series as $E$. We here introduce the standard shorthand. Given an integer $q \geq 1$ and a central character $\chi \pmod q$, denote by $Y_0(q) = \Gamma_0(q)\backslash\mathbb{H}$ the modular curve and by $L^2(Y_0(q), \chi)$ the space of measurable functions $f : \mathbb{H} \to \mathbb{C}$ satisfying the periodicity condition

$$f\left(\frac{az+b}{cz+d}\right) = \chi(d)f(z) \quad \text{for all} \quad \begin{pmatrix} a & b \\ c & d \end{pmatrix} \in \Gamma_0(q)$$

and $\langle f, f\rangle_q < \infty$, where $\langle \cdot, \cdot\rangle_q$ signifies the Petersson inner product

$$\langle f, g\rangle_q := \int_{Y_0(q)} f(z)\overline{g(z)}d\mu(z)$$

with the Poincaré measure $d\mu(z) = y^{-2}dxdy$. We also indicate by $\mathcal{A}(Y_0(q), \chi)$ the space of smooth automorphic functions of nebentypus $\chi$. For a prime number $q$, the spectrum of $L^2(Y_0(q))$ consists of the constant function, Maaß forms, and unitary Eisenstein series $E_\infty(\cdot, 1/2 + it)$, $E_0(\cdot, 1/2 + it)$ for $t \in \mathbb{R}$, corresponding to the two ($\Gamma_0(q)$-equivalence classes of) cusps $\mathfrak{a} \sim \infty, 0$. Of our two cusps in question, one is termed *narrow* and the other *wide*.

The chief hindrance to reprove Theorem 1.3 arises out of the fact that $|E|^2 \notin L^2(Y_0(q))$. By assembling our workarounds described below in detail, the correct version of Theorem 1.3 is stated as follows:

**Theorem 1.3.** *Let $q$ traverse the set of primes and let $t \in \mathbb{R}^\times$. With the notation above, for a fixed $\phi$ of level 1, we then have*

$$\int_{Y_0(q)} \phi(z)|E_\infty(z, 1/2 + it)|^2 d\mu(z) \longrightarrow \frac{6\log q}{\pi}\int_{Y_0(1)} \phi(z)d\mu(z).$$

We refer the reader to the recent paper of Pan–Young [12] for the treatment of more general context, as our limit law in Theorem 1.3 is deduced from [12, Theorem 1.2] with $q = 1$. As shall be seen in Remark 1.1, the level $q$ Eisenstein series associated to the wide cusp cannot satisfy any reasonable form of QUE.

*Proof.* We sketch out the overall strategy of the proof. Our present reasoning completely differs from the original argument of the second author, since the usage of incomplete Eisenstein series necessitates a huge load of computations. We go through the regularisation procedure of Zagier [17] (or its adèlic treatment by Michel–Venkatesh [9]). We make do with a similar analysis to [12] based upon the construction of $\mathcal{E} \in \mathcal{A}(Y_0(q))$ such that $|E|^2 - \mathcal{E} \in L^2(Y_0(q))$. Here $\mathcal{E}$ is a linear combination of Eisenstein series of level $q$ and trivial central character. Let $\mathcal{H}_{it}(q, 1)$ stand for the set of $L^2(Y_0(q))$-normalised Hecke–Maaß newforms of level $q$, trivial central character, and spectral parameter $t$. Via the spectral decomposition roughly of the form

$$\langle |E|^2 - \mathcal{E}, \phi\rangle_q \approx \sum_{u_j \in \mathcal{H}_{it_j}(1,1)} \sum_{t_j \ll T} \langle |E|^2, u_j\rangle_q \langle u_j, \phi\rangle_1 + (\text{continuous spectrum}),$$

we calculate as follows:

$$\begin{aligned}\langle |E|^2, \phi\rangle_q &= \langle |E|^2 - \mathcal{E}, \phi\rangle_q + \langle \mathcal{E}, \phi\rangle_q \\ &= \lim_{\eta\to 0^+} \Big( \langle E_\infty(\cdot, 1+\eta), \phi\rangle_q + \varphi_{\infty 0}(1/2+it)\varphi_{\infty 0}(1/2-it+\eta) \\ &\qquad \langle E_0(\cdot, 1-\eta), \phi\rangle_q \Big) + O(\cdots),\end{aligned}$$

where $O(\cdots)$ signifies an inexplicit but acceptable error term and $\varphi$ is an entry of the scattering matrix. We centre upon the resulting quantity of the limit as $\eta \to 0$ with emphasis on the fact that $E_\infty(z, 1+\eta)$ and $E_0(z, 1-\eta)$ have simple poles at $\eta = 0$ with respective residues $\pm\mathrm{vol}(Y_0(q))^{-1}$ independently of $z \in Y_0(q)$. By virtue of unitarity of the scattering matrix (cf. [5, (13.30)]), the poles cancel with each other. An exact calculation of the scattering matrix yields

$$\varphi_{\infty 0}(1/2+it)\varphi'_{\infty 0}(1/2-it) = -2\log q + O((\log\log(q+1))^5);$$

thence the constant term in the above limit turns into

$$\frac{2\langle 1, \phi\rangle_q}{\mathrm{vol}(Y_0(q))}\log q + \langle \Gamma_0(q)\text{-automorphic function}, \phi\rangle_q.$$

It suffices for our purpose to deal only with the main term. We are able to connect $\langle\cdot,\cdot\rangle_q$ with an integration on level 1 by folding the integrals, obtaining

$$\langle |E|^2, \phi\rangle_q = \frac{2\nu(q)\langle 1, \phi\rangle_1}{\mathrm{vol}(Y_0(q))}\log q + \langle \mathrm{SL}_2(\mathbb{Z})\text{-automorphic function}, \phi\rangle_1 + O(\cdots), \quad (1.1)$$

where $\nu(q) = [\mathrm{SL}_2(\mathbb{Z}), \Gamma_0(q)]$. We may specify the magnitude of the resulting error term to deduce a more accurate asymptotic formula. Hence, recalling $\nu(q) = q\prod_{p|q}(1+p^{-1}) = 3\mathrm{vol}(Y_0(q))/\pi$ and that the second term on the right-hand side of (1.1) can be small enough, we are led to what we wanted to prove. We do not go into details, for the rigorous proof of Theorem 1.3 and the generalisation to all Eisenstein newforms of squarefree level or of primitive central character are elucidated in [12]. □

Having proven Theorem 1.3, we observe that Theorem 1.2 in the original paper [7] should hold (for the narrow cusp $\mathfrak{a} \sim \infty$), although its proof was erroneous. By taking $\phi$ to be an indicator function of the set $A_1$ in Theorem 1.3 formulated above, the asymptotic formulae derived in these two theorems will have the same left-hand side up to measures.

In the aforementioned proof, we have relied on Zagier's theory of renormalisation. We could, however, tailor the Luo–Sarnak method to our setting to reobtain Theorem 1.3, because the spectral decomposition of $\phi$ is pursued on level 1 and there are no forms of level $q$ that enter. This fixes any problem of trying to adapt the Luo–Sarnak method for varying levels. As to the germinal methodology in [7, pp. 1134–1149], we came to the conclusion that the error terms in Propositions 1.4 and 1.5 are not provably smaller than the main terms; it was forgotten to ponder over the approximation of an indicator function of $A_q$.

For squarefree levels, *every* Eisenstein series is a linear combination of oldforms. We contend that QUE should centre upon newforms (or alternatively $p$-adic microlocal lifts in the language of Nelson [11]). Instead of picking the usual orthonormal basis $\{E_{\mathfrak{a}}(z, 1/2+it) : \mathfrak{a} \text{ is a cusp of } Y_0(q)\}$, we concoct an orthonormal basis out of Eisenstein newforms and oldforms (see [16, Sect. 8.4]). In this sense, the Eisenstein newforms $E_{\chi_1,\chi_2}$ would be the sheer counterpart of the holomorphic case, and Theorem 1.3 above might not be the most felicitous analogue of [6]. Here $\chi_i$ is primitive $\pmod{q_i}$, $i = 1, 2$ and $q = q_1q_2$. As a proper generalisation of [7], Pan–Young [12] studied Eisenstein series for general nebentypen and levels with $(q_1, q_2) = 1$, and used the explicit computation of the Fourier expansion of Eisenstein series and Zagier's theory of renormalisation. Specifically, their results imply that QUE turns out to be ruined for all varying test functions (designed to fulfil some conventions) at least if $M \gg q^{\delta}$ for some $\delta > 0$, where $M$ is the level of the test function with $M \mid q$. One then perceives that the statement in [12, Sect. 1.5] contradicts [7, Theorems 1.2 and 1.3]. Nonetheless, our Theorem 1.3 amounts to a stepping stone towards solving the *genuine* QUE problem in the level aspect.

*Remark 1.1.* In hindsight, QUE for the wide cusp $\mathfrak{a} \sim 0$ could not be true in view of the change-of-basis formula of Young [16, Theorem 6.1] or Booker–Lee–Strömbergsson [2]. To be more precise, one obtains

$$E_0^{(q)}(z, s) = q^{-s}\frac{\zeta(2s)}{L(2s, \chi_{0,q})}\{E(z, s) - q^{-s}E(qz, s)\}, \tag{1.2}$$

where $E$ stands for the classical Eisenstein series of level 1. Putting $\Re(s) = 1/2$, we see that $E_0$ differs mildly as $q \to \infty$ (by the relative error of $q^{-1/2}$ in the Eisenstein norm) from the corresponding level 1 Eisenstein series. Theorems 1.2 and 1.3 in [7] therefore give vacuous statements for the wide cusp, because $E(z, s)$ in (1.2) is independent of $q$.

## 2. Modifications of Main Computation

Next we point out problematic errors on pp. 1139–1143, and correct the argument imitating the notion of Luo–Sarmak [8] at the sacrifice of a certain restriction on $\phi$. The following modifications are meant for the original paper [7] to have a mostly self-contained exposition. First of all, the definition of the Fourier coefficients $\rho_j(n)$ on p. 1139 is in error; it should be replaced with

$$\rho_j(n) := \begin{cases} q^{1/2}\tau_j(n/q) & \text{if } q \mid n, \\ 0 & \text{otherwise,} \end{cases} \tag{2.1}$$

and this propagates typos throughout the later parts of the paper when we treat oldforms of level 1 slashed by $q$. For adherent miscalculations coming from this change, we refer to the itemised list in Sect. 3. In particular, the assertion (2) in Lemma 2.4 is incorrect and we reform it to the shape that includes the extra $q^{1/2}$.

**Lemma 2.4.** (2) *Let $u_j$ be an oldform for $\Gamma_0(q)$ expressed as $u_j(z) = (\iota_q v_j)(z)$ with $v_j$ a cusp form for* $\mathrm{SL}_2(\mathbb{Z})$*, and $(\iota_\ell f)(z) := f(\ell z)$ for every $\ell \geq 1$. Then we have*

$$\sum_{n=1}^{\infty}\frac{\rho_j(n)\sigma_\nu(n)}{n^s} = \tau_j(1)q^{1/2-s}\frac{1 + q^\nu - \tilde{\tau}_j(q)q^{\nu-s}}{1 - q^{\nu-2s}}\frac{L(s, v_j)L(s-\nu, v_j)}{\zeta(2s-\nu)}$$

*and*

$$\sum_{n=1}^{\infty}\frac{\rho_j(n)\sigma_\nu(nq^{-\alpha})}{n^s} = \tau_j(1)q^{1/2-s}(1 + q^\nu)\frac{1 - \tilde{\tau}_j(q)q^{\nu-s}}{1 - q^{\nu-2s}}\frac{L(s, v_j)L(s-\nu, v_j)}{\zeta(2s-\nu)}.$$

Here we have also simplified the complicated numerator in (2.13) as well. A thing to be aware of is that the factor $L(s, \nu_j)L(s-\nu, \nu_j)\zeta(2s-\nu)^{-1}$ is akin to the Rankin–Selberg $L$-function handled in [12]:

$$\sum_{n=1}^{\infty} \frac{\overline{\lambda_{1,\overline{\psi}}}(n)\lambda_j(n)}{n^{1/2+it}} = \frac{L(1/2, u_j)L(1/2+2it, u_j \otimes \overline{\psi})}{L(1+2it, \overline{\psi} \otimes \chi_{0,q})}.$$

Here $\lambda_{\chi_1,\chi_2}(n) = \chi_2(\operatorname{sgn}(n)) \sum_{ab=|n|} \chi_1(a)\overline{\chi_2}(b)(b/a)^{it}$ and $\lambda_j(n)$ is the $n$th Fourier coefficient of the Hecke–Maaß newform $u_j$ of level $q$ and spectral parameter $t_j$. In view of the altered version of Lemma 2.4 (2), it would be better to combine Lemma 2.6 with Lemma 2.5 to conclude the calculation of the contribution from cusp forms, since the proof of Lemma 2.5 contains many incorrect statements (Lemma 2.6 is back-of-the-envelope). For example, on p. 1143, the second-to-last line in the fourth display should be modified into

$$\begin{aligned}
&\frac{\rho_j(1)L(s+it, u_j)L(s-it, u_j)}{(1-q^{-2s})\zeta(2s)\hat{\zeta}(1+2it)}\left(1 - \frac{1-\tilde{\rho}_j(q)q^{-s-it}}{1-q^{-1-2it}}\right)\\
&= \frac{\rho_j(1)L(s+it, u_j)L(s-it, u_j)}{(1-q^{-2s})\zeta(2s)\hat{\zeta}(1+2it)}\Big(1-(1-\tilde{\rho}_j(q)q^{-s-it})(1+q^{-1-2it}+\cdots)\Big)\\
&\ll \frac{\rho_j(q)q^{-s-it}}{1-q^{-2s}} \frac{L(s+it, u_j)L(s-it, u_j)}{\zeta(2s)\hat{\zeta}(1+2it)} \ll q^{17\vartheta/15-1/10+\epsilon} \ll q^{0.02396},
\end{aligned}$$

where we exploited the subconvex bound in the work of Blomer et al. [1]: $L(1/2, u_j) \ll q^{1/5+\vartheta/15+\epsilon}$. We denote by $\vartheta \in [0, 7/64]$ an admissible exponent towards the Ramanujan–Petersson conjecture and note that the appearance of $\vartheta$ stems also from $\rho_j(q)$. However, the sum over $j$ is missing in [7] and we must consider the second moment of $L(1/2+it, u_j)$. This process will render a power saving error term. The Stade formula (Lemma 3.1) coupled with the Stirling asymptotics implies that the integral in the last display on p. 1142 decays exponentially in $t$ with $s = 1/2+it$, namely

$$\int_0^{\infty} K_{it}(2\pi y)K_{it_j}(2\pi y)y^{s-1}dy \ll \mathcal{P}(t, t_j) \times \begin{cases} 1 & \text{if } |t_j| \leq |t|/2, \\ \exp(\pi(|t|/2-|t_j|)) & \text{if } |t_j| > |t|/2, \end{cases}$$

for a certain polynomial $\mathcal{P}$ in $t$ and the spectral parameter $t_j$. Then we use the Lindelöf-on-average bound

$$\sum_{t_j \leq |t|/2+\log q} \sum_{u_j \in \mathcal{H}_{it_j}(1,1)} |L(1/2+it, u_j)|^2 \ll_{t,\epsilon} q^{\epsilon},$$

which follows readily from the spectral large sieve. Here we implicitly supposed that $\phi$ is a test function on level 1, and thus the hindrance in evaluating the contribution of cusp forms was that we assumed $\phi$ to be of level $q$.

In addition, the full part of the last big display (2.15) on p. 1143 is mistaken. In order to rectify it, we use Lemma 2.4 (2), obtaining the following plausible estimate:

$$\begin{aligned}(2.15) &= \frac{1}{\hat{\zeta}(1+2it)}\Bigg(\sum_{n=1}^{\infty}\frac{\rho_j(n)\sigma_{-2it}(n)}{n^{s-it}} - \frac{1}{1-q^{-1-2it}}\sum_{n=1}^{\infty}\frac{\rho_j(n)\sigma_{-2it}(nq^{-\alpha})}{n^{s-it}}\Bigg)\\ &= \tau_j(1)\frac{q^{1/2+it-s}}{1-q^{-2s}}\frac{L(s+it,v_j)L(s-it,v_j)}{\zeta(2s)\hat{\zeta}(1+2it)}\\ &\quad\times\Big(1+q^{-2it}-\tilde{\tau}_j(q)q^{-s-it}-(1+q^{-2it})(1-\tilde{\tau}_j(q)q^{-s-it})(1+q^{-1-2it}+\cdots)\Big)\\ &= \tau_j(q)\frac{q^{1/2-2(s+it)}}{1-q^{-2s}}\frac{L(s+it,v_j)L(s-it,v_j)}{\zeta(2s)\hat{\zeta}(1+2it)} + \text{(lower order terms)}\\ &\ll q^{\vartheta-1/2}.\end{aligned}$$

Hence, the contribution of oldforms would not cause any problem in analysing the integral $I_j(q,s)$. In any case, the original reduction to the asymptotic evaluation of $\langle |E|^2, E_{\mathfrak{a}}(\cdot|\psi)\rangle_q$ for the incomplete Eisenstein series $E_{\mathfrak{a}}(z|\psi)$ cannot work at all unless $\phi$ is a test function on $Y_0(1)$.

## 3. List of Other Mistakes

The following list provides local corrections to numerous mathematical, notational, and peripheral errors. These modifications allow us to reasonably restate rudimentary materials prepared in [7]. The notation here is completely the same as in the original paper and may differ from that in the foregoing discussions.

(1) In the introduction, we incorrectly used the phrase of quantum ergodicity. In this regard, results in [3,13,18] state that on a closed Riemannian manifold whose geodesic flow is ergodic with respect to the Liouville measure, a density one subsequence of eigenfunctions converges microlocally to this measure. We wanted to generally mean QUE in phase space or configuration space.

(2) p. 1132: the coefficient on the right-hand side of the first display should be $6/\pi$ in the light of Stade's formula on $\mathrm{GL}_2$; see [14, (1.1)] for a correct statement. That is, when we consider the modular surface equipped with the Poincaré measure $d\mu(z)$, one then has

$$\int_{Y_0(1)}\phi(z)|E(z,1/2+it)|^2 d\mu(z) = \frac{3}{\pi}\log(1/4+t^2)\int_{Y_0(1)}\phi(z)d\mu(z) + o(\log t)$$

as $t\to\infty$ for a *fixed* smooth and compactly supported test function $\phi : Y_0(1)\to\mathbb{R}$. In passing, Young [15] has considered the rate of convergence with power saving error terms, obtaining the following sophisticated formula:

$$\langle |E|^2,\phi\rangle = \frac{\langle 1,\phi\rangle_1}{\langle 1,1\rangle_1}\log(1/4+t^2) + O\Big(\frac{\|\phi\|_1\log t}{\log\log t}\Big) + O(\|\phi\|_2\, t^{-1/6+\epsilon}).$$

(3) p. 1132 and thereafter: the cusps $\infty$ and $i\infty$ represent the same and they are to be unified.

(4) pp. 1132–1134, Theorems 1.2, 1.3 and Propositions 1.4, 1.5 are all incorrect. These errors necessitate thorough corrections, as described below.

(5) p. 1133, the second display: the volume $\mathrm{vol}(A^0)$ ought to be written as $\mathrm{vol}(A^0_q)$, although the asymptotic formula stated no longer holds.
(6) p. 1134, Proposition 1.5: we adopted the notation $V_q(z)$ and $V_q$ to mean the same object, and they must be unified to either. The same happens in Proposition 3.1.
(7) p. 1134, (2.1): the integral in the first line should be over $\Gamma_{q,\kappa}\backslash H^2$, though the method does not work for the wide cusp.
(8) p. 1135 and thereafter: it should be stated that $\kappa$ is short for $\kappa(q)$ in order to avoid any confusion.
(9) p. 1135, (2.3) and (2.4): the condition on the sum should be recasted as $c \in \mathcal{C}(\kappa, \kappa')$, where $\mathcal{C}(\kappa, \kappa')$ is a set of allowed moduli attached to a pair of cusps $(\kappa, \kappa')$ and associated scaling matrices $\sigma_\kappa, \sigma_{\kappa'}$:

$$\mathcal{C}(\kappa, \kappa') = \left\{ \gamma > 0 : \begin{pmatrix} * & * \\ \gamma & * \end{pmatrix} \in \sigma_\kappa^{-1}\Gamma_0(q)\sigma_{\kappa'} \right\}.$$

In particular, we have ([5, Sect. 4.2])

$$\begin{aligned} \mathcal{C}(\infty, \infty) = \mathcal{C}(0, 0) &= \{\gamma = \ell q : \ell \in \mathbb{N}\}, \\ \mathcal{C}(\infty, 0) = \mathcal{C}(0, \infty) &= \{\gamma = c\sqrt{q} : c \in \mathbb{N}, (c, q) = 1\}. \end{aligned} \tag{3.1}$$

We redefine the Kloosterman sum to a modulus $\gamma \in \mathcal{C}(\mathfrak{a}, \mathfrak{b})$ for arbitrary cusps $\mathfrak{a}, \mathfrak{b}$ by

$$S_{\mathfrak{a}\mathfrak{b}}(m, n; \gamma) := \sum_{\left(\begin{smallmatrix} a & b \\ \gamma & d \end{smallmatrix}\right) \in \Gamma_\infty \backslash \sigma_{\mathfrak{a}}^{-1}\Gamma_0(q)\sigma_{\mathfrak{b}} / \Gamma_\infty} e\left(\frac{am + dn}{\gamma}\right).$$

When $\mathfrak{a} = \mathfrak{b} \sim \infty$ we obtain the usual Kloosterman sum $S(m, n; c)$ for $c \equiv 0 \pmod{q}$ with dropping the subscripts. For the pair of cusps $(\infty, 0)$, the formula (3.1) gives $S_{\infty 0}(m, n; c\sqrt{q}) = S(\overline{q}m, n; c)$. In this regard, the sum of Kloosterman sums turns for $q$ prime into

$$\sum_{c \in \mathcal{C}(\kappa, \kappa')} \frac{S_{\kappa\kappa'}(0, n; c)}{c^{2s}} = \begin{cases} \displaystyle\sum_{(c,q)=1} \frac{S(0, n; c)}{q^s c^{2s}} & \text{for the } (\infty, 0) \text{ cusp-pair}, \\ \displaystyle\sum_{q|c>0} \frac{S(0, n; c)}{c^{2s}} & \text{for the } (\infty, \infty) \text{ cusp-pair}, \end{cases}$$

with $S(0, n; c)$ equal to the Ramanujan sum $\sum_{d|(n,c)} \mu(c/d)d$. For our study of Eisenstein series, we are largely interested in the case of $\kappa = \kappa'$; we thus have either the $(\infty, \infty)$ or $(0, 0)$ cusp-pair. On the other hand, via an explicit computation of scaling matrices, we infer that $\sigma_\infty^{-1}\Gamma_0(q)\sigma_\infty = \sigma_0^{-1}\Gamma_0(q)\sigma_0 = \Gamma_0(q)$. By the above reformulation, the change of variables $c = \gamma\sqrt{q}$ on p. 1138 is justified.
(10) p. 1135, (2.4): for the additional factor $n^{s-1}$ read $|n|^{s-1}$ with $n \in \mathbb{Z} \setminus \{0\}$. Correspondingly, the same correction also applies to $n^{s-1}$ in the last display on p. 1137, and to $n^{-1/2+it}$ in the first display on p. 1143.
(11) p. 1135: the notation for the Kloosterman sum is not unified around its definition; namely $\mathcal{S}_{q,\kappa,\kappa'}$ and $S_{q,\kappa,\kappa'}$ indicate the same function.

(12) p. 1135, (2.6): the condition $1 \le d < c$ will not suffice and it should be $d \pmod c$, or $1 \le d \le c$. Accordingly, the first sum in the second display on p. 1136 should be

$$\sum_{d \,(\mathrm{mod}\, c)}^{*} e\Big(\frac{nd}{c}\Big).$$

Notice that we have also changed $e(na/c)$ to $e(nd/c)$. Likewise the condition $1 \le \delta < \gamma$ in the fifth display on p. 1136 should be replaced with $\delta \pmod \gamma$.

(13) p. 1136, line 2: for $\mathcal{S}_{q,\kappa,\infty}(0,n,c)$ read $\mathcal{S}_{q,\kappa,\kappa'}(0,n,c)$ in the original notation.

(14) p. 1136, line 7: for $B\backslash\sigma_a^{-1}\Gamma_0(q)\sigma_{\kappa'}/B$ read $B\backslash\sigma_\kappa^{-1}\Gamma_0(q)\sigma_{\kappa'}/B$.

(15) p. 1136, (2.7): it is sufficient for our requirements to write $\begin{pmatrix}\alpha\sqrt{q} & * \\ \gamma\sqrt{q} & \delta\sqrt{q}\end{pmatrix}$ instead of $\begin{pmatrix}\alpha\sqrt{q} & \beta\sqrt{q} \\ \gamma\sqrt{q} & \delta\sqrt{q}\end{pmatrix}$ because we do not impose any restriction on $\beta$. In order to include $\beta$, the upper-right entry is to be of the form $\beta/\sqrt{q}$ so that $\alpha\delta q - \beta\gamma = 1$ with $\alpha, \beta, \gamma, \delta \in \mathbb{Z}$.

(16) p. 1136, the third and fifth displays: we do not need to write $\in \sigma_0^{-1}\Gamma_0(q)$ and this is in fact incorrect. It should also be remarked that $\sigma_\infty^{-1}\Gamma_0(q)\sigma_0 = \sigma_0^{-1}\Gamma_0(q)\sigma_\infty$.

(17) p. 1136, Lemma 2.2: here is a rather easier proof of this claim. For notational convenience, we reexpress the Ramanujan sum as $r_c(n) = a_n(c)$, then

$$\sum_{c=1}^{\infty}\frac{r_c(n)}{c^s} = \sum_{c=1}^{\infty}\sum_{\substack{d_1|n\\ d_1d_2=c}}\frac{d_1\mu(d_2)}{(d_1d_2)^s} = \sum_{d_2=1}^{\infty}\frac{\mu(d_2)}{d_2^s}\sum_{d_1|n}d_1^{1-s} = \frac{\sigma_{1-s}(n)}{\zeta(s)}.$$

This finishes the proof.

(18) p. 1137, Lemma 2.3: the quantity $\alpha_p$ should read $\alpha_q$ and the $c$-sum should be over $c \in \mathcal{C}(\kappa,\kappa')$. In addition we deal with $n \neq 0$, thereby we can redesign Lemma 2.3 to be valid for all $n \neq 0$.

(19) p. 1138, (2.9): for the $K$-Bessel function $K_{ir_j}(2\pi ny)$ read $K_{ir_j}(2\pi|n|y)$, and similarly $K_{ir_j}(2\pi ny)$ and $K_{ir_j}(2\pi qny)$ on p. 1139 should be $K_{ir_j}(2\pi|n|y)$ and $K_{ir_j}(2\pi q|n|y)$, respectively. The same should apply to $K_{it}(2\pi ny)$ in the second last and the last displays on p. 1146.

(20) p. 1139, the fourth display: it should be that, with our taking the reformulation (2.1) into account,

$$\sum_{n=1}^{\infty}\frac{\rho_j(n)}{n^s} = \sum_{n=1}^{\infty}\frac{\rho_j(qn)}{(qn)^s} = q^{1/2-s}\tau_j(1)L(s,v_j).$$

In the same way, the factor $q^{1/2}$ should enter into the first and second displays on p. 1142.

(21) p. 1141, the last display: in the second line, the first term $1/(1-\tilde{\tau}(q)q^s)$ is to be recasted as $1/(1-\tilde{\tau}(q)q^{-s})$. Correspondingly, the same correction applies to the third line in the first display on p. 1142. The shape of the resulting $L$-functions is nevertheless correct.

(22) p. 1142, the last display: we have the Mellin transform of the two $K$-Bessel functions and now state the basic lemma to simplify the claim in Lemma 2.5.

**Lemma 3.1** (Stade's formula for $\mathrm{GL}_2$ [4, (6.576.4)]). *The Mellin transform of the product of $K$-Bessel functions can be evaluated as follows:*

$$\int_0^\infty K_\mu(y)K_\nu(y)y^s\frac{dy}{y} = 2^{s-3}\Gamma\Big(\frac{s+\mu+\nu}{2}\Big)\Gamma\Big(\frac{s+\mu-\nu}{2}\Big)$$
$$\times\Gamma\Big(\frac{s-\mu+\nu}{2}\Big)\Gamma\Big(\frac{s-\mu-\nu}{2}\Big)\Gamma(s)^{-1}.$$

We find the computation in [8] was erroneous in the sense that the additional factor in front of the gamma functions did not appear in their paper (we make the change of variables $y \mapsto 2\pi y$ to apply Lemma 3.1).

(23) p. 1144, line 2: the estimate $\tilde{\tau}_j(q) = O(q^{1/4})$ was employed, while we generally have $\tilde{\tau}_j(q) = O(q^\vartheta)$.

(24) p. 1144, the proof of Proposition 2.6: we exploited the path-breaking subconvex bound due to Duke–Friedlander–Iwaniec. We notice that their bound has nowadays been improved in the subsequent literature, especially in connection with the Motohashi formula [10].

(25) p. 1145, the fourth display: it is required to insert the extra factor of 2 in front of the sum over $n \geq 1$, as the Fourier expansion (2.2) of Eisenstein series entails the sum over $n \in \mathbb{Z} \setminus \{0\}$. The same error propagates through the argument in the sequel; specifically the last display on p. 1146 must read

$$I_2 = \frac{8}{|\hat{\zeta}(1+2it)|^2}\int_{(2)} H(s)\sum_{n=1}^\infty \frac{|\eta_{q,\kappa,\kappa}(n,1/2+it)|^2}{n^s}\int_0^\infty |K_{it}(2\pi y)|^2 y^{s-1}\frac{dsdy}{2\pi i},$$

where $|K_{it}(2\pi ny)|^2$ has been altered into $|K_{it}(2\pi y)|^2$.

(26) p. 1145, the second last display: the quantity $2\Re(\varphi_{q,\kappa,\kappa}(1/2+it)y^{2it})$ in the numerator ought to be $2\Re(\varphi_{q,\kappa,\kappa}(1/2+it)y^{-2it})$. The same correction may apply to the fourth display on p. 1146.

(27) p. 1146, the first display: it would be desirable to write

$$\sum_{q|c}\frac{\varphi(c)}{c^{2w}} = \frac{\zeta(2w-1)}{\zeta(2w)}\Big(1-\frac{1-q^{1-2w}}{1-q^{-2w}}\Big) = \frac{\zeta(2w-1)}{\zeta(2w)}\frac{q-1}{q^{2w}-1}.$$

(28) p. 1147, (3.1): the function $|\eta_{q,\kappa,\infty}(n,1/2+it)|^2$ on the left-hand side should read $|\eta_{q,\kappa,\kappa}(n,1/2+it)|^2$. The same should be pursued for the first display on p. 1148. Besides, the coefficient in front of the $\alpha$-sum on the right-hand side in (3.1) should be

$$\frac{\zeta(s)^2\zeta(s+2it)\zeta(s-2it)(1-q^{-s})(1-q^{2it-s})(1-q^{-2it-s})}{\zeta(2s)(1+q^{-s})|1-q^{-2it}|^2|1-q^{-1-2it}|^2}.$$

(29) p. 1147, the last display: the resulting shape of the sum is indeed

$$4\sum_{\alpha=0}^\infty \frac{\sin^2\theta}{q^{\alpha s}} = \frac{q^{-s}(1+q^{-s})|1-q^{-2it}|^2}{(1-q^{2it-s})(1-q^{-2it-s})(1-q^{-s})},$$

where $\theta = t\alpha\log q$. With this expression in mind, the first display on p. 1148 becomes

$$\sum_{n=1}^\infty \frac{|\eta_{q,\kappa,\kappa}(n,1/2+it)|^2}{n^s} = q^{-s}\frac{\zeta(s)^2\zeta(s+2it)\zeta(s-2it)}{\zeta(2s)|1-q^{-1-2it}|^2} + O(q^{-1}).$$

The same process will be applied to the second and second last displays on p. 1148.

(30) p. 1148: we have to remedy the argument as well as minor mistakes. We initially invoke Lemma 3.1 above, reconstructing $I_2$ as

$$I_2 = \frac{1}{|\hat{\zeta}(1+2it)|^2}\int_{(2)} B(s)\frac{ds}{2\pi i} + O(q^{-1}),$$

where we put

$$B(s) = H(s)\frac{\zeta(s)^2\zeta(s+2it)\zeta(s-2it)}{\zeta(2s)|1-q^{-1-2it}|^2}\frac{\Gamma(s/2)^2\Gamma(s/2+it)\Gamma(s/2-it)}{(\pi q)^s\Gamma(s)}.$$

We then shift the contour to $\Re(s) = 1/2$ to derive

$$I_2 = \frac{\mathrm{res}_{s=1}B(s)}{|\hat{\zeta}(1+2it)|^2} + \frac{1}{|\hat{\zeta}(1+2it)|^2}\int_{(1/2)} B(s)\frac{ds}{2\pi i} + O(q^{-1}).$$

This is justified by the Stirling formula and the fact that $H(\sigma+it)$ is of rapid decay in $t$. The integral over $\Re(s) = 1/2$ converges and produces $O(1)$ as $q \to \infty$. Hence, we arrive at

$$I_2 = \frac{\mathrm{res}_{s=1}B(s)}{|\hat{\zeta}(1+2it)|^2} + O(1).$$

(31) p. 1148, the last display: the error term ought to be $O(|s-1|)$.

(32) p. 1149, the third display: there must be a minus sign in front of $\log q$.

(33) p. 1149, the second last display: the normalised measure was defined on p. 1131 as $dV = \mathrm{vol}(X)^{-1}y^{-2}dxdy$. One can therefore use the Poincaré measure $d\mu(z)$ to denote $y^{-2}dxdy$.

(34) In References:
- [DFI]: “artin” ⟶ “Artin” and also add “(with an appendix by Dorian Goldfeld, Jeffrey Hoffstein, and Daniel Lieman)” to the title.
- [HL]: “maass” ⟶ “Maass” and “siegel” ⟶ “Siegel”.
- [I]: add the year of publication “2002”.
- [K]: “eisenstein” ⟶ “Eisenstein”.
- [KMV]: “Rankin-selberg” ⟶ “Rankin–Selberg”.
- [RS]: “behavior” ⟶ “behaviour”.
- [PS]: $SL(2, o)\backslash H^3$ ⟶ $\mathrm{SL}_2(\mathcal{O})\backslash\mathbf{H}^3$.
- [S]: “rankin-selberg” ⟶ “Rankin–Selberg”; besides “419–445” ⟶ “419–453”.

*Acknowledgements.* The authors are much indebted to Jiakun Pan for discussions and Matthew Young for informing them of various slips in the original paper. Thanks are also owed to Peter Humphries for his feedback and the anonymous referees for suggestions that helped improve this article.